\documentclass[12pt]{article}%

\usepackage{amssymb}
\usepackage{amsfonts}
\usepackage{amsmath}
\usepackage[nohead]{geometry}
\usepackage[singlespacing]{setspace}
\usepackage[bottom]{footmisc}
\usepackage{indentfirst}
\usepackage{endnotes}
\usepackage{graphicx}%
\usepackage{rotating}
\usepackage{verbatim}
\usepackage{setspace}
\usepackage{multirow}
\usepackage{latexsym}
\usepackage{subfigure}

\usepackage{amsthm}
\usepackage{makeidx}
\usepackage{fancyhdr}
\usepackage{type1cm}

\newcommand{\hV}{\widehat  V }

\newcommand{\tf}{\widetilde f}

\newcommand{\tF}{\widetilde F}

\newcommand{\hC}{\widehat C}
\newcommand{\hF}{\widehat F}
\newcommand{\hf}{\widehat f}

\newcommand{\hu}{\widehat u}

\newcommand{\hSig}{\widehat\Sigma}
\newcommand{\tSig}{\widetilde\Sigma}
\newcommand{\hsig}{\widehat\sigma}
\newcommand{\hlam}{\widehat\lambda}
\newcommand{\hLam}{\widehat \Lambda}

\newcommand{\cov}{\mathrm{cov}}

\newcommand{\diag}{\mathrm{diag}}

\newcommand{\var}{\mathrm{var}}
\newcommand{\beq}{\begin{eqnarray*}}
\newcommand{\eeq}{\end{eqnarray*}}

\newcommand{\Sigu}{\Sigma_{u}}

\numberwithin{equation}{section}
\theoremstyle{plain}
\newtheorem{thm}{Theorem}[section]

\newtheorem{lem}{Lemma}[section]

\newtheorem{prop}{Proposition}[section]
\newtheorem{assum}{Assumption}[section]
\theoremstyle{definition}
\newtheorem{exm}{Example}[section]
\newtheorem{remark}{Remark}[section]

\makeatletter
\def\@biblabel#1{\hspace*{-\labelsep}}
\makeatother
\begin{document}

\title{Statistical Inferences Using Large Estimated Covariances for Panel Data and Factor Models }
\author{Jushan Bai\footnote{Department of Economics, Columbia University, New York, NY 10027.   }\medskip\\{Columbia University}  \and Yuan Liao\footnote{Department of Mathematics, University of Maryland at College Park, College Park, MD 20742.}\medskip\\{University of Maryland}   }
\date{}
\maketitle

\sloppy%

\onehalfspacing


\begin{abstract}

While most of the convergence results in the literature on high dimensional
covariance matrix are concerned about the accuracy of
estimating the covariance matrix (and precision matrix), relatively
less is known about the effect of estimating large covariances on statistical
inferences. We study two important models: factor analysis
and panel data model with interactive effects, and focus on the statistical
inference and estimation efficiency of structural parameters
based on   large covariance estimators.   For efficient estimation, both models call for
a weighted principle components (WPC), which relies on a high dimensional
weight matrix. This paper derives an efficient and
feasible WPC using the covariance matrix estimator of Fan et al. (2013).
  However,  we demonstrate that existing results on large covariance  estimation based on absolute convergence are not suitable for statistical  inferences of the structural parameters. What is needed is some weighted consistency and the associated rate of convergence, which are obtained in this paper.  Finally, the proposed
method is applied to the US divorce rate data. We find that the efficient WPC
identifies the significant effects of divorce-law reforms on the divorce
rate, and it provides more accurate estimation and tighter confidence intervals than  existing methods.

\end{abstract}
\strut

\strut

\pagebreak%
\doublespacing

\onehalfspacing

\section{Introduction}


Estimating a high-dimensional covariance matrix has been an active research area in the recent literature.  Many methods are proposed for estimating the covariance matrix and the precision (inverse covariance) matrix, e.g.  El Karoui (2008), Bickel and Levina (2008), Rothman et al. (2009),  Lam and Fan (2009), Cai and Liu (2011), Fan et al. (2013). Among many theoretical results, rates of convergence under various interesting matrix norms have been derived. In particular, if we write $N$ to denote the dimension and $T$ to denote the sample size, when the $N\times N$ covariance matrix $\Sigma$ is sparse whose  eigenvalues are bounded away from zero,  we can  obtain an estimator $\hSig$ that achieves  a \textit{near-$\sqrt{T}$}-rate under the \textit{operator norm}:
\begin{equation}\label{e1.1}\|\hSig-\Sigma\|=O_p(m_N(\frac{\log N}{T})^{\frac{1-q}{2}})=\|\hSig^{-1}-\Sigma^{-1}\|
\end{equation}
where $m_N$ and $q$ are parameters that measure the level of sparsity. Cai and Zhou (2012)  showed that the rate of convergence (\ref{e1.1}) is minimax optimal. However, there is relatively less knowledge about the effect of estimating a high-dimensional covariance matrix on statistical inferences, e.g., the estimation efficiency for a parametric model, and the effect of estimating large covariances on   the limiting distributions for estimators of some structural parameters.

We find  that  when a high-dimensional covariance estimator is applied  for  statistical inferences (precisely, deriving limiting distributions of estimated structural parameters), most of the results in the literature based on \textit{absolute convergence} like (\ref{e1.1}) are not suitable, even with the  minimax optimal   rate. Instead, a   ``\textit{weighted convergence}" is needed,   which takes the form $\|A_1(\hSig^{-1}-\Sigma^{-1})A_2\|$, where both $A_1, A_2$ are stochastic matrices that weight the estimation error $\hSig^{-1}-\Sigma^{-1}$. The weights $A_1$ and $A_2$ further ``average down" the estimation errors, which significantly improve the rate of convergence    to make valid statistical inferences.  However, the weighted  convergence cannot be implied by the usual results   in the   literature.  One of our contributions is to tackle this problem.

This paper focuses on two    models that are of increasing importance in  many statistical applications: factor analysis and panel data model with interactive effects.   In  factor analysis,  the notion of sparsity     is a natural assumption based on the factor structure, which is proved to be a successful approach   (e.g.,  Boivin and Ng 2006, Phan 2012, Andersen et al. 2011).  This paper  gives a theoretical justification about how such a sparse structure can be used  to improve the estimation efficiency in two general models.   Both problems involve estimating a large weight matrix, where the problem of proving ``weighted convergence" is present.


\subsection{Approximate factor model}
We consider a high-dimensional approximate factor model:
\begin{equation}\label{e12}
y_{it}=\lambda_{i}'f_t+u_{it},\quad i\leq N, t\leq T.
\end{equation}
where $f_t$ is an $r\times 1$ vector of  common factors, $\lambda_i$ is a vector of factor loadings, and $u_{it}$ represents the error term, often known as the \textit{idiosyncratic component}.  If we denote $Y_t=(y_{1t},...,y_{Nt})'$, $\Lambda=(\lambda_{1},....,\lambda_{N})'$, and $u_t=(u_{1t},...,u_{Nt})'$,  model (\ref{e12}) can  be written as
$$
Y_t=\Lambda f_t+u_t,\quad t\leq T.
$$ Only $Y_{t}$ is observable in the model.  In a data-rich environment, both $N$ and $T$ can be large and the dimension $N$ might be even much larger than $T$. The goal is to make efficient inference about  $\lambda_i'f_t$, $\lambda_i, f_t$ or their rotations.


Approximate factor models often require the $N\times N$ covariance matrix $\Sigma_u=\cov(u_t)$ be non-diagonal matrix and the diagonal entries may vary over a large range (Chamberlain and Rothschild 1983).  The traditional 
  method of principal components  (PC)    essentially treats $u_{it}$ to be  homoskedastic and uncorrelated over $i$.  As a result,   it is    inefficient. 
In this paper, we consider a \textit{weighted principal components } (WPC) method to efficiently estimate the   heteroskedastic approximate factor models. The WPC   solves a  weighted least squares problem:
  \begin{equation}\label{e1.2}
\min_{\Lambda, f_t}\sum_{t=1}^T(Y_t-\Lambda f_t)'W(Y_t-\Lambda f_t)
\end{equation}
subject to certain normalization constraints.   Here $W$ is an $N\times N$ positive definite weight matrix.  We propose  a feasible  efficient WPC that  requires  consistently estimating the high-dimensional $\Sigma_u^{-1}$ (when $N>T$) as the weight matrix, and  is shown to be optimal over a broad class of estimators.

\subsection{Large panel data model with interactive effects}
A closely related model  is the panel data   with a factor structure in the error term:
\begin{equation}\label{e1.3}
y_{it}=X_{it}'\beta+\varepsilon_{it},\quad \varepsilon_{it}=\lambda_i'f_t+u_{it},\quad i\leq N, t\leq T,
\end{equation}
where $X_{it}$ is a $d\times 1$ vector of regressors; $\beta$ is a $d\times 1$ vector of unknown coefficients. The regression noise   $\varepsilon_{it}$ has a factor structure with unknown loadings and factors.  In the model, the only observables are $(y_{it}, X_{it})$. The goal is to estimate the structural parameter  $\beta$,  whose dimension is fixed. In this model, the factor component $\lambda_i'f_t$ is regarded as an \textit{interactive effect} of   the individual and time effects.  Because the regressor and factor can be correlated,   simply regressing $y_{it}$ on $X_{it}$ is not consistent.

Similarly, we propose  to  estimate $\beta$ via:   \begin{equation}\label{e1.4}
\min_{\beta, \Lambda, f_t}\sum_{t=1}^T(Y_t-X_{it}'\beta-\Lambda f_t)'W(Y_t-X_{it}'\beta-\Lambda f_t),
\end{equation}
with a high-dimensional weight matrix $W$.   The method is also WPC because the  estimated factors are shown to be principal components of the weighted sample covariance matrix.   In particular, it allows a consistent estimator for $\Sigma_u^{-1}$ as the optimal weight matrix even when $\Sigma_u^{-1}$ is non-diagonal and $N/T\rightarrow\infty$. 
Except for sparsity, the  off-diagonal structure of $\Sigma_u$ is unknown. The WPC takes into account both cross-sectional correlation and heteroskedasticity of $u_{it}$ over $i$, while the existing methods in the literature, e.g., 
  Bai 2009, Moon and Weidner 2010,  do not.

 \subsection{Summary of contributions}


 First of all,    we develop the inferential theory  using a general high-dimensional weight $W$. This admits many promising choices of the weight matrices that are suitable for specific applied problems for factor analysis. Especially, in cases where estimating $\Sigma_u$ is difficult, our inferential theory  is still useful when suitable weight matrices are chosen to improve the estimation efficiency.   Secondly, we show that when $W=\Sigma_u^{-1}$  is used, the WPC yields an optimal estimator  in the sense that the estimated common component  $\lambda_i'f_t$ and structural parameter $\beta$ have the minimum  asymptotic variance over a broad class of estimators. 
 
 Third,  we focus on the effect of estimating large covariance matrices on  efficient statistical inferences.  
In both pure factor analysis  and the large panel data with  a factor structure,   we employ a consistent estimator for $\Sigma_u^{-1}$ recently proposed by Fan et al. (2013), as an operational  weight matrix.   Therefore, our optimal estimator is  still feasible under $N/T\rightarrow\infty$.  However, substituting a consistent estimator $\Sigma_u^{-1}$   is highly non-trivial when $N>T$. An interesting phenomenon is observed:   most existing results on estimating large covariances are not suitable for statistical inferences of the models being considered.    We develop a new strategy that investigates the \textit{weighted consistency} for the estimated optimal weight matrix to address this problem.
 


Fourth, we consistently estimate the  asymptotic variances of the  proposed estimators   under both cross-sectional and serial correlations in $u_{it}$. Hence the new WPC estimator for the interative effect model is readily  used for statistical inferences in practice.  In contrast,   existing methods usually require additionally modeling the large error covariance (e.g., assuming diagonality, parametrizing the off-diagonal structure) in order for practical inferences.

Given  the popularity of the PC method,  why do we need a new estimator to incorporate the large covariance $\Sigma_u$?
Even though most of the existing methods for panel data models avoid estimating $\Sigma_u$, to demonstrate the potential   efficiency loss for existing methods, we present a real-data application in Section \ref{law}, which studies the effect of divorce reform law on the change of  divorce rates. The WPC is applied to the  year-state divorce rate data of U.S. during 1956-1985. It illustrates that after incorporating $\Sigma_u^{-1}$ in the estimation, WPC captures the significant (negative) effects from nine to twelve years after the law was reformed, consistent with the previous empirical findings in the social science literature. In contrast, the existing method (PC) without estimating $\Sigma_u^{-1}$ would result in wide confidence intervals and potentially conservative conclusions.  Numerically, we find an average of  46\%  efficiency gained using WPC, relative to the existing method.  In addition, the proposed WPC also enjoys the    computational convenience, as it also admits analytical solutions.


 Realizing the limitation of the regular PC method, some  important works have been developed   to improve the estimation efficiency for factor analysis, e.g., Breitung and Tenhofen (2011),   Bai and Li (2012) and Doz et al. (2012). They   require the cross-sectional dependences' structure  be  specifically modeled.   Recently, Choi (2012) specified $W=\Sigma_u^{-1}$, which  essentially  requires $\Sigma_u$ be known.  
 Recently, Fan et al. (2013) proposed a thresholding method   to estimate $\Sigma_u^{-1}$. They focused on covariance matrix estimations and   did not address the efficient estimation for the factors, loadings  and panel data models. As we   discussed, replacing $\Sigma_u^{-1}$ with its consistent estimator is technically challenging when $N/T\rightarrow\infty$.   Additional  literature on factor analysis and panel data with interactive effects includes, e.g.,  Pesaran (2006),  Ahn et al. (2001), Su and Chen (2013),  Su et al. (2012),   Wang (2009),  Forni et al. (2000),     Hallin and Li\v{s}ka (2007), Lam and Yao (2012),  Cheng and Hansen (2013),  Caner and Han (2012),  etc. None of these  incorporated $\Sigma_u^{-1}$ or studied  efficient estimation for panel data models.   We also remark that there is a rapidly growing literature on estimating high-dimensional (inverse) covariance matrices. Besides those mentioned, the list also includes, e.g.,  d' Aspremont et al. (2008), Bien and Tibshirani (2011), Luo (2011),   Pati et al. (2012), Xue et al. (2012),  among many others.

We assume  the number of factors $r=\dim(f_t)$ to be known. When $r$ is unknown, it can be consistently estimated by certain information criteria as in, e.g., Bai and Ng (2002), as we shall briefly discuss in Section \ref{secIn}. 

 \vspace{2em}



 The rest of the paper is organized as follows. Section 2 describes the general problem of statistical inference based on large covariance matrices. Section 3  formally proposes the  WPC method. The large-sample inferential theory  of WPC with a general weight matrix  is presented.   Section \ref{sepc}  introduces the  efficient WPC.   Section \ref{secIn} applies the WPC method to the panel data model with interactive effects.   Section \ref{ssim} illustrates numerical  comparisons of related methods. Section \ref{law} applies WPC to a real data problem of divorce rate study. Finally, Section 8 concludes. All proofs are given in the supplementary material.

Throughout the paper, we use $\lambda_{\min}(A)$ and $\lambda_{\max}(A)$ to denote the minimum and maximum eigenvalues of matrix $A$. We also let $\|A\|$, $\|A\|_1$ and $\|A\|_F$  denote the operator norm, $L_1$-norm and Frobenius norm of a matrix, defined as $\|A\|=\sqrt{\lambda_{\max}(A'A)}$, $\|A\|_1=\max_{i}\sum_{j}|A_{ij}|$ and $\|A\|_F=\sqrt{\sum_{i,j}A_{ij}^2}$ respectively. Note that if $A$ is a vector, $\|A\|=\|A\|_F$ is equal to the Euclidean norm. Finally, for two sequences, we write  $a_T\ll b_T$ (and equivalently $b_T\gg a_T$) if $a_T=o(b_T)$ as $T\rightarrow\infty.$

\section{Challenge of Inference based on Large Estimated Covariance}

 
  
Consider estimating a low-dimensional structural parameter $\theta$ that arises from a   model involving a high-dimensional covariance matrix $\Sigma$. It is often the case that when $\Sigma$ were known, incorporating it in the estimator may achieve a better estimation accuracy, e.g., smaller standard errors and  tighter confidence intervals. Taking into account $\Sigma$, the estimator can be written as a function of the data $D_T$ and   $\Sigma$ as ($T$ denotes the sample size):$$\widehat\theta=f(D_T,\Sigma),$$ and the limiting distribution may be derived. In practice, we replace $\Sigma$   by a consistent estimator $\hSig$ and obtain  a feasible efficient estimator  $f(D_T, \hSig)$.

To show that replacing $\Sigma$ with its consistent estimator does not affect the limiting distribution of $\widehat\theta$, one often needs 
$a_T(f(D_T,\Sigma)-f(D_T,\hSig))=o_p(1)$ where $a_T^{-1}$ can be understood as the rate of convergence of $\widehat\theta$.  However, such a simple substitution is  technically difficult if $N>T$. To see this, note that often  $f(D_T,\Sigma)$ depends on the precision matrix $\Sigma^{-1}$, and  the effect of estimating $\Sigma^{-1}$ is approximately linearly dependent on $\hSig^{-1}-\Sigma^{-1}$. We can often write 
$$
a_T(f(D_T,\Sigma)-f(D_T,\hSig))=a_TA_1(\hSig^{-1}-\Sigma^{-1})A_2+o_p(1)
$$where $A_1, A_2$ are typically non-sparse stochastic matrices  of dimensions $\dim(\theta)\times N$ and $N\times 1$ respectively.     Applying the  Cauchy-Schwarz inequality,   $$
a_T\|A_1(\hSig^{-1}-\Sigma^{-1})A_2\|\leq a_T\|A_1\|\|A_2\|\|\hSig^{-1}-\Sigma^{-1}\|.
$$
As both $A_1$ and $A_2$ are high-dimensional   matrices (vectors), the right hand side of the above inequality is typically not  stochastically negligible even if  the ``absolute convergence" $\|\hSig^{-1}-\Sigma^{-1}\|$ achieves the optimal  convergence rate.\footnote{When $\Sigma$ is sparse enough, one can obtain a near $\sqrt{T}$-rate of convergence for the $L_1$-norm $\|\hSig^{-1}-\Sigma^{-1}\|_1$, but this still yields a crude bound for $a_TA_1(\hSig^{-1}-\Sigma^{-1})A_2$.} The problem arises because  $\|A_1\|$ and $\|A_2\|$ grow fast with the dimensionality, so they accumulate the estimation errors and lead to a crude bound.

We further illustrate this issue in two examples, which  are to be studied in detail in  the paper. 

\begin{exm}
Consider the high-dimensional factor model (\ref{e12}). The parameter of interest is the common component $\lambda_i'f_t$. The efficient estimation crucially depends on $
\frac{1}{\sqrt{N}}\Lambda'(\hSig_u^{-1}-\Sigma_u^{-1})u_t,
$
for a sparse covariance estimator $\hSig_u^{-1}$.  However, the existing  results on the optimal convergence of $\|\hSig_u^{-1}-\Sigma_u^{-1}\|$ in the literature (e.g., Fan et al. 2013) are not applicable  directly when $N>T$,  because  $\|\Lambda\|=O(\sqrt{N})$ and $\|u_t\|=O_p(\sqrt{N})$, but the minimax rate for $\|\hSig_u^{-1}-\Sigma_u^{-1}\|$  is no faster than $O_p(T^{-1/2})$. 
Applying the absolute convergence for $\widehat\Sigma_u^{-1}$, 
$\frac{1}{\sqrt{N}}\|\Lambda\|\|\widehat\Sigma_u^{-1}-\Sigma_u^{-1}\|\|u_t\|=O_p(\sqrt{\frac{N}{T}})\neq o_p(1)$ when $N>T$.
\end{exm}

\begin{exm}
Consider the high-dimensional panel data model (\ref{e1.3}).    The efficient estimation of $\beta$  requires estimating the inverse covariance $\Sigma_u^{-1}$. Suppose   $\widetilde\Sigma_u^{-1}$ is a consistent estimator. We require
$$
\frac{1}{\sqrt{NT}}Z'[(\tSig_u^{-1}-\Sigma_u^{-1})\otimes I_T]U=o_p(1),
$$
where $I_T$ is a $T$-dimensional identity matrix and $Z$ and $U$ are stochastic matrices whose dimensions are $NT\times \dim(\beta)$ and $NT\times 1$ respectively. However, because $\|Z\|=O_p(\sqrt{NT})=\|U\|$, it is difficult to apply the  absolute convergence $\|\tSig_u^{-1}-\Sigma_u^{-1}\|$ (whose minimax rate is no faster than $O_p(T^{-1/2})$)  to achieve the desired convergence when  $N>T$. 
The crude bound gives $\frac{1}{\sqrt{NT}}\|Z\|\|\widetilde\Sigma_u^{-1}-\Sigma_u^{-1}\|\|U\|=O_p(\sqrt{N})\neq o_p(1)$. $\square$
\end{exm}

As one of the main contributions of this paper,   a new strategy of   ``weighted convergence" is developed. When analyzing  $a_TA_1(\widehat\Sigma^{-1}-\Sigma^{-1})A_2$, we should not separate the covariance estimation error from the weighting matrices $A_1, A_2.$ Intuitively, the weights further ``average down" the estimation errors, to ensure the asymptotic negligibility of   the weighted error.  We demonstrate that the weighted  convergence is useful for high-dimensional inferences in panel data and factor models, and cannot be simply implied by the usual results on ``absolute convergence" in the   literature.      
  
  
\section{Approximate  Factor Models}

\subsection{Weighted principal components}

In model (\ref{e12}), the only observables are $\{Y_t\}_{t=1}^T$, and both the factors $\{f_t\}_{t=1}^T$ and loadings $\Lambda=(\lambda_1,...,\lambda_N)'$ are parameters to  estimate.   
 We estimate them via  the following weighted least squares:
\begin{equation}\label{eqq2.1}
(\widehat{\Lambda}, \hf_t)=\min_{\Lambda, f_t}\sum_{t=1}^T(Y_t-\Lambda f_t)'W_T(Y_t-\Lambda f_t)
\end{equation}
 subject to:
\begin{equation}\label{eq2.1}
\frac{1}{T}\sum_{t=1}^T\hf_t\hf_t'=I_r;\quad\hLam'W_T\hLam\text{ is diagonal}.
\end{equation}
Here $W_T$ is an $N\times N$ weight matrix, which can be either stochastic or deterministic. When $W_T$ is stochastic, we   mean $W_T$ to be  a consistent estimator of some   positive definite  $W$ under the operator norm.     
We will show in Section \ref{sepc} that the optimal weight   is $\Sigma_u^{-1}$.   On the other hand, keeping a general $W_T$ admits other choices of  the weight for specific applied problems, especially when it is difficult to estimate the optimal weight matrix.



Solving (\ref{eqq2.1}) subjected to  the restriction  (\ref{eq2.1}) gives the WPC estimators:    $\hlam_j$ and $\hf_t$ are both ${r}\times 1$ vectors such that,   the columns of the $T\times {r}$ matrix $\hF/\sqrt{T}=(\hf_1,...,\hf_T)'/\sqrt{T}$ are the eigenvectors  corresponding to the largest ${r}$ eigenvalues of $Y'W_TY$, and $\hLam=T^{-1}Y\hF=(\hlam_1,...,\hlam_N)'.$ We call our method to be \textit{weighted principal components} (WPC), to distinguish from   the traditional principal components (PC) method  that uses $W_T=I_N$. Note that  PC does not encounter the problem of estimating large covariance matrices, and is not efficient when $\{u_{it}\}$'s are cross-sectionally correlated across $i$.

It has been well known that  the factors and loadings are not separably identifiable without further restrictions. 
The WPC   estimates rotated factors and loadings with  rotation matrix $H_W$. Let $\hV$ be the $r\times r$ diagonal matrix of the first $r$ largest eigenvalues of $YW_TY'/(TN)$. Let $F=(f_1,...,f_T)'$, then  $H_W=\hV^{-1}\hF'F\Lambda'W_T\Lambda/(NT)$. We use the subscript $W$ to emphasize the dependence of the rotation on $W$.


\subsection{General conditions}\label{scr}
We present general results for the proposed WPC with a general weight matrix, which hold for a broad class of estimators.    
For the general weight matrix $W$ and its data-dependent version $W_T$,   the following assumption is needed:

 \begin{assum}\label{a26} (i) $\|W_T-W\|=o_p(\min\{T^{-1/4}, N^{-1/4},  \sqrt{\frac{N}{T}}, \sqrt{\frac{T}{N\log N}}\})$.\\
 (ii) $\|\frac{1}{\sqrt{N}}\Lambda'(W_T-W)u_t\|=o_p(1).$
 \end{assum}

Condition (i) is easy to satisfy by using many ``good" covariance estimators given in the literature. 
However,    the main challenge described in Section 2 arises from proving   condition (ii) in the above assumption. When $W_T$ is a consistent estimator for $\Sigma_u^{-1}$, we shall see in Section \ref{swc} that  this requires a new ``weighted convergence", which is   necessary but challenging to the high-dimensional inference problems being considered.

  We allow the factors and idiosyncratic components to be weakly serially dependent  via the strong mixing condition.   Let $\mathcal{F}_{-\infty}^0$ and $\mathcal{F}_{T}^{\infty}$ denote the $\sigma$-algebras generated by $\{(f_t,u_t): -\infty\leq t\leq 0\}$ and  $\{(f_t,u_t): T\leq t\leq \infty\}$ respectively. In addition, define the mixing coefficient
\begin{equation} \label{mixing}
\alpha(T)=\sup_{A\in\mathcal{F}_{-\infty}^0, B\in\mathcal{F}_{T}^{\infty}}|P(A)P(B)-P(AB)|.
\end{equation}

\begin{assum}\label{a21} (i)
 $\{u_t, f_t\}_{t\geq1}$ is strictly stationary. In addition, $Eu_{it}=Eu_{it}f_{jt}=0$ for all $i\leq p, j\leq r$ and $t\leq T.$
\\
(ii) There exist  constants $c_1, c_2>0$   such that  $c_2<\lambda_{\min}(\Sigma_u)\leq\lambda_{\max}(\Sigma_u)<c_1,$  $\max_{j\leq N}\|\lambda_{j}\|<c_1$, and $c_2<\lambda_{\min}(\cov(f_t))\leq\lambda_{\max}(\cov(f_t))<c_1$.
\\
(iii)   Exponential tail:   There exist $r_1, r_2>0$ and $b_1, b_2>0$, such that for any $s>0$, $i\leq p$ and $j\leq r$,
$P(|u_{it}|>s)\leq\exp(-(s/b_1)^{r_1}), $ and $ P(|f_{jt}|>s)\leq \exp(-(s/b_2)^{r_2}).$\\
(iv) Strong mixing: There exists   $r_3>0$  and $C>0$ satisfying:  for all $T\in\mathbb{Z}^+$,
$$\alpha(T)\leq \exp(-CT^{r_3}).$$
\end{assum}

 We  assume that $W$ has bounded row sums, that is, $\|W\|_1<M$ for some $M>0$. Write $\Lambda'W=(d_1,...,d_N)$, with each $d_i$ being an $r\times 1$ vector.  Then  $\max_{j\leq N}\|d_j\|<\infty.$ 
 
   The following assumptions    are standard in the literature.
Assumption \ref{a22}   requires the factors  be \textit{pervasive}, which  holds when the factors  impact a non-vanishing proportion of individual time series.  Assumption  \ref{a27} extends similar conditions in Stock and Watson (2002) and Bai (2003). When  $W=I_N$ is used, they reduce to those in the literature of the regular PC. A simple sufficient condition for Assumption \ref{a27} is that $u_{it}$ is i.i.d. in both $i$ and $t$.






 \begin{assum}\label{a22} All the eigenvalues of the  $r\times r$  matrix $\Lambda'\Lambda/N$ are  bounded away from both zero and infinity as $N\rightarrow\infty$.
\end{assum}

\begin{assum}\label{a27}
(i)  $E\|\frac{1}{\sqrt{NT}}\sum_{s=1}^Tf_s(u_s'Wu_t-Eu_s'Wu_t)\|^2=O(1)$.\\
(ii) For each $i\leq N$, $E\|\frac{1}{\sqrt{NT}}\sum_{t=1}^T\sum_{j=1}^Nd_j(u_{jt}u_{it}-Eu_{jt}u_{it})\|=O(1)$.\\
(iii) For each $i\leq r$, $E\|\frac{1}{\sqrt{NT}}\sum_{t=1}^T\sum_{j=1}^Nd_ju_{jt}f_{it}\|=O(1)$.\\
(iv) There is a constant $\delta\geq 4$ and $M>0$ such that for all large $N$,\\$E|\frac{1}{\sqrt{N}}(u_s'Wu_t-Eu_s'Wu_t)|^{\delta}<M$ and  $E\|\frac{1}{\sqrt{N}}\Lambda'Wu_t\|^{\delta}<M$.
\end{assum}


 \subsection{Limiting distributions}
 
 
 

The factors and loadings are two sets of parameters to estimate. The limiting distributions of their estimators depend on the following asymptotic expansions, to be shown   in the appendix: for some positive definite matrix   $J_W$,  and the rotation matrix $H_W$, 
\begin{eqnarray}\label{eq2.8}
&&\sqrt{N}(\hf_t-H_Wf_t)=J_W\frac{\Lambda'Wu_t}{\sqrt{N}}+O_p(a_{T})\cr 
&&\sqrt{T}(\hlam_j-H_W^{'-1}\lambda_j)=H_W\frac{1}{\sqrt{T}}\sum_{t=1}^Tf_tu_{jt}+O_p(b_{T}).
\end{eqnarray}
where the asymptotic normality arises from the leading terms while $a_{T}$ and $b_{T}$ are some remaining stochastic sequences.


The limiting distribution of $\hlam_j$ requires  $H_W$ to have a limit. We thus need the following condition:

  \begin{assum} \label{a28}
  (i) There is an $r\times r$ matrix $\Sigma_{\Lambda}$ such that $\Lambda'W\Lambda/N\rightarrow\Sigma_{\Lambda}$ as $N\rightarrow\infty$. In addition,  the eigenvalues of the  $\Sigma_{\Lambda}\cov(f_t)$ are distinct.  
 \\
 (ii) For each $t\leq T$,
$(\Lambda'W\Sigma_uW\Lambda)^{-1/2}\Lambda'Wu_t\rightarrow^d\mathcal{N}(0,I_r).$
  \end{assum}

   According to the expansions of (\ref{eq2.8}),   the above condition (ii)     is almost a necessary condition for the asymptotic normality of $\widehat f_t$.   Note that  $\frac{1}{\sqrt{N}}\Lambda'Wu_t=\frac{1}{\sqrt{N}}\sum_{i=1}^Nd_iu_{it}$. Hence a cross-sectional central limit theorem can indeed apply.      Condition (ii) is only for $\widehat f_t$,  and  the  limiting distribution of the estimated   loading $\widehat\lambda_j$ in Theorem \ref{t33} below  does not depend on this condition.




  We now introduce some notation that are needed to present the limiting distributions. Let $V$ be an $r\times r$ diagonal matrix with element as the largest $r$ eigenvalues of $\Sigma_{\Lambda}^{1/2}\cov(f_t)\Sigma_{\Lambda}^{1/2}$, and $\Gamma_W$ be the corresponding eigenvector matrix such that $\Gamma_W'\Gamma_W=I_r$. We use the subscript $W$ to indicate that $\Gamma_W$ depends on $W$ via $\Sigma_{\Lambda}$. Recall that $\Sigma_{\Lambda}$ is defined in Assumption \ref{a28}.  Let $Q_W=V^{1/2}\Gamma_W'\Sigma_{\Lambda}^{-1/2}$. In fact $H_W\rightarrow^p Q_W^{'-1}$. In addition, to account for the serial correlation over $t$,  let 
\begin{equation}\label{eq2.9}
\Phi_j=E(f_tf_t'u_{jt}^2)+\sum_{t=1}^{\infty}E[(f_1f_{1+t}'+f_{1+t}f_1')u_{j1}u_{j,1+t}].
\end{equation}

\begin{thm} \label{t33} Assume  $(\log N)^2=o(T)$ and $T=o(N^2)$. Under Assumptions \ref{a26}-\ref{a28}(i), for each $j\leq N$, 
$$
\sqrt{T}(\hlam_j-H_W^{'-1}\lambda_j)\rightarrow^d \mathcal{N}(0, Q_W^{'-1}\Phi_jQ_W^{-1}).
$$
If in addition,   $N=o(T^2)$ and Assumption  \ref{a28}(ii) holds, 
$$
N(V^{-1}Q_W\Lambda'W\Sigma_uW\Lambda Q_W'V^{-1})^{-1/2}(\hf_t-H_Wf_t)\rightarrow^d \mathcal{N}(0,I_r).
$$
For the common component, we have $$
\frac{\hlam_i'\hf_t-\lambda_i'f_t}{(\lambda_i'\Xi_W\lambda_i/N+f_t'\Omega_if_t/T)^{1/2}}\rightarrow^d\mathcal{N}(0,1).
$$
where $\Xi_W=\Sigma_{\Lambda}^{-1}\Lambda'W\Sigma_uW\Lambda\Sigma_{\Lambda}^{-1}/N
$ and  $\Omega_i=\cov(f_t)^{-1}\Phi_i\cov(f_t)^{-1}$. 
\end{thm}

  \begin{remark}
The eigenvalues of $(V^{-1}Q_W\Lambda'W\Sigma_uW\Lambda Q_W'V^{-1})^{-1/2}$ are of order $O(N^{-1/2})$.  Hence Theorem \ref{t33} implies the $\sqrt{N}$-consistency of the estimated factors. If  we further assume that $\Lambda'W\Sigma_uW\Lambda/N $ has a limit, say  $G$,  then immediately we have
$$
\sqrt{N}(\hf_t-H_Wf_t)\rightarrow^d \mathcal{N}(0, V^{-1}Q_WGQ_W'V^{-1}),
$$
where the $\sqrt{N}$-consistency is more clearly demonstrated.
\end{remark}

The uniform convergence of $\hf_t$ and $\hlam_j$ are given below.

\begin{thm} \label{t21} Let  $\alpha=\max\{1/r_1, 1/r_2\}$ with $r_1, r_2$ defined in Assumption \ref{a21}. Let $\delta\geq4$ be as defined in Assumption  \ref{a27}.  Under  Assumptions \ref{a26}-\ref{a27},  as $N, T\rightarrow\infty$,
\begin{eqnarray}
 \label{eq2.5}
&&\max_{t\leq T}\|\hf_t-H_Wf_t\|=O_p\left((\log T)^{\alpha}\|W_T-W\|+\frac{T^{1/\delta}}{\sqrt{N}}+\frac{1}{\sqrt{T}}\right),\\
&&\label{eq2.6}\max_{j\leq N}\|\hlam_j-H_W^{'-1}\lambda_j\|=O_p\left(\|W_T-W\|+\frac{1}{\sqrt{N}}+\sqrt{\frac{\log N}{T}}\right).
\end{eqnarray}

\end{thm}

 

 \begin{remark}
  

The uniform convergence in (\ref{eq2.5}) and (\ref{eq2.6}) is important under large $N$ and $T$.  For example, in estimating large covariance matrices, it is used to derive the proper levels of  thresholding  or shrinkage (e.g., Fan et al. 2013, Ledoit and Wolf 2012).

 \end{remark}





\subsection{Heteroskedastic WPC}

As a simple choice for $W$, 
$$W= (\diag(\Sigma_u))^{-1}.$$ 
This choice  improves the regular PC when cross-sectional heteroskedasticity is present. This weight can be easily estimated  using the residuals. First apply the regular PC by taking $W_T=I_N$, and obtain a  consistent estimator $\hC_{it}$ of the common component $\lambda_i'f_t$ for each $i\leq N, t\leq T.$ Define $$
W_T^h=\diag\{\hsig_{u,11}^{-1},...,\hsig_{u,NN}^{-1}\},\text{ where } \hsig_{u,ii}=\frac{1}{T}\sum_{t=1}^T(y_{it}-\hC_{it})^2.
$$
Then in the second step, apply the WPC  with   weight matrix $W_T^h.$  

The heteroskedastic WPC (which we call HWPC) method has been previously suggested by, e.g.,  Breitung and Tenhofen (2011). Investigations of its theoretical properties  can be found in the appendix. Moreover, numerical studies in Section \ref{ssim}  show that this method improves the efficiency relative to the regular PC method.

 \section{Efficient Principal Components Under Conditional Sparsity}\label{sepc}


In the approximate factor models,  $u_{it}$'s are   correlated (over $i$).   A more efficient estimator (which we call EWPC)  should take $W=\Sigma_u^{-1}$ as the weight matrix. This estimator has been recently suggested by Choi (2012), but   $\Sigma_u^{-1}$ was assumed to be known. 

There are two main challenges in practice: on one hand, when $N>T$, $\Sigma_u^{-1}$ is hard to estimate as the sample covariance based on the residual $\widehat u_{it}$ is no longer invertible.  On the other hand, as we illustrated in Section 2, even if a consistent estimator for $\Sigma_u^{-1}$ is available, it is technically difficult to prove that the effect of covariance estimation is neglibile when $N/T\rightarrow\infty.$ We first     apply Fan et al. (2013)'s method to estimate $\Sigma_u^{-1}$, and then  address the second problem in Section \ref{swc}.

\subsection{Conditional Sparsity}\label{scs}

We apply a thresholded covariance estimator to estimate  $\Sigma_u^{-1}$, which is recently proposed by   Fan et al. (2013) for factor analysis. 
 Let $(\nu_j,\xi_j)_{j=1}^N$ be the eigenvalues-eignvectors of the sample covariance $S_y$ of $Y_t$, in a decreasing order such that $\nu_1\geq \nu_2\geq...\geq\nu_N.$  Let 
$$
R=S_y-\sum_{i=1}^r\nu_i\xi_i\xi_i'.
$$
Define a general thresholding function $s_{ij}(z): \mathbb{R}\rightarrow\mathbb{R}$ as in Rothman et al. (2009) and Cai and Liu (2011)  with an entry-dependent threshold $\tau_{ij}$ such that:\\
(i) $s_{ij}(z)=0$  if $|z|<\tau_{ij};$\\
(ii) $|s_{ij}(z)-z|\leq \tau_{ij}.$\\
(iii) There are constants $a>0$ and $b>1$ such that $|s_{ij}(z)-z|\leq a\tau_{ij}^2$ if $|z|>b\tau_{ij}$.\\
Examples of $s_{ij}(z)$ include the hard-thresholding: $s_{ij}(z)=zI_{(|z|>\tau_{ij})}$; SCAD (Fan and Li 2001), MPC (Zhang 2010) etc.   As for the threshold value, we specify
\begin{equation}\label{e42}
\tau_{ij}= C\sqrt{R_{ii}R_{jj}}\omega_T,\text{ where } \omega_T= \sqrt{\frac{\log N}{T}}+\frac{1}{\sqrt{N}}
\end{equation}
for some pre-determined universal $C>0$, chosen from cross-validation as in  Fan et al. (2013).  Then estimate $\Sigma_u$ by $\hSig_u=(\hSig_{u,ij})_{N\times N}$,
$$
\hSig_{u,ij}=\begin{cases}
R_{ii}, & i=j\\
s_{ij}(R_{ij}), & i\neq j
\end{cases}, \quad\text{ where } R=(R_{ij})_{N\times N}.
$$
Intuitively, $\hSig_u$ thresholds off the small entries of the residual covariance $\frac{1}{T}\sum_{t=1}^T\hu_t\hu_t'$ obtained from the regular PC estimate.

 To apply such a weight estimator,  we assume $\Sigma_u$ to be a sparse matrix. In an approximate factor model,      such a special structure is known to be \textit{conditionally sparse} (given the common factors). 
 Consider the notion of \text{generalized sparsity}: write $\Sigma_{u}=(\Sigma_{u,ij})_{N\times N}$. For some $q\in[0,1/2)$, define
\begin{equation}\label{e40}
m_N=\max_{i\leq N}\sum_{j=1}^N|\Sigma_{u,ij}|^q.
\end{equation}
In particular, when $q=0$, define  $m_N=\max_{i\leq N}\sum_{j=1}^NI_{(\Sigma_{u,ij}\neq0)}$. Mathematically, the conditional sparse structure on $\Sigma_u$ assumes, there is   $q\in [0,1/2)$, such that
\begin{equation}\label{e41}
m_N=o\left(\min\left\{\frac{1}{N^{1/4}}\left(\frac{T}{\log N}\right)^{(1-q)/2}, N^{1/4-q/2}\right\}\right).
\end{equation}

In the sparse covariance estimation literature, Condition (\ref{e41}) itself is enough to achieve a covariance estimator such that $\|\Sigma_u^{-1}-\hSig_u^{-1}\|=o_p(1)$, whose rate of convergence is nearly $\sqrt{T}$ (e.g.,    Cai and Zhou 2012, Fan et al. 2013, etc.). But for the ``weighted convergence" needed for efficient estimations in factor analysis and large panel data models, this condition is not sufficient. Therefore,  we  introduce a more refined description of the sparse structure of $\Sigma_u$ (condition (ii) in Assumption \ref{a41} below), which is similar to those in Rothman  et al. (2008). 
 
 Let $S_L$ and $S_U$ denote two disjoint sets and  respectively include the indices of small and large elements of $\Sigu$ in absolute value, and
\begin{equation}\label{eq4.3}
\{(i,j): i\leq N, j\leq N\}=S_L\cup S_U.
\end{equation}
We assume $(i,i)\in S_U$ for all $i\leq N.$  The    sparsity condition  assumes that most of the indices $(i,j)$ belong to $S_L$ when $i\neq j$.  A special case arises when $\Sigma_{u}$ is strictly sparse, in the sense that its elements with small magnitudes ($S_L$) are exactly zero. For the banded matrix as an example,  $ \Sigma_{u,ij}=0 $ if $|i-j|> k$
for some fixed $k\geq 1.$ Then $S_L=\{(i,j):|i-j|>k\}$ and $S_U=\{(i,j):|i-j|\leq k\}$.  Another example is the block-diagonal matrix.

The following assumption mathematically defines the ``conditional sparsity"   for the approximate factor model.

Define  $$\omega_T=\sqrt{\frac{\log N}{T}}+\frac{1}{\sqrt{N}}.$$
\begin{assum}\label{a41}
(i) There is $q\in[0,1/2)$ such that (\ref{e41}) holds.\\
(ii) There is a partition $\{(i,j): i\leq N, j\leq N\}=S_L\cup S_U$ such that  $\sum_{i\neq j, (i,j)\in S_U}1=O(N)$ and  
  $\sum_{(i,j)\in S_L}|\Sigma_{u,ij}|=O(1)$. In addition, 
$$\max_{(i,j)\in S_L}|\Sigma_{u,ij}|=O(\omega_T), \quad \omega_T=O(\min_{(i,j)\in S_U}|\Sigma_{u,ij}|).$$
\end{assum}

If for example, $\Sigma_u$ is a block covariance matrix with finite block sizes, this assumption is naturally satisfied as long as the signal is not too-weak (that is,  $\omega_T=o(\min_{(i,j)\in S_U}|\Sigma_{u, ij}|)$).  
     Condition (ii)  requires the elements in $S_L$ and $S_U$ be well-separable.  The partition $\{(i,j): i\leq N, j\leq N\}=S_L\cup S_U$ may not be unique.  Most importantly, we do not need to know either $S_L$ or $S_U$; hence the block size, the banding length, or the locations of the zero entries can be completely unknown.     Our analysis suffices as long as such a partition exists.

\subsection{Weighted convergence using the optimal weight matrix}\label{swc}

We now formally discuss the issue brought by Assumption \ref{a26}. In order for the effect of estimating $\Sigma_u^{-1}$ to be negligible,   $\|\frac{1}{\sqrt{N}}\Lambda'(\hSig_u^{-1}-\Sigma_u^{-1})u_t\|=o_p(1)$ is required, which is a tight condition. However, a  direct application of the optimal rate of convergence (i.e., Fan et al. 2013, Cai and Zhou 2012) $\|\hSig_u^{-1}-\Sigu^{-1}\|=O_p(m_N\omega_T^{1-q})$ implies 
 $$
 \|\frac{1}{\sqrt{N}}\Lambda'(\hSig_u^{-1}-\Sigma_u^{-1})u_t\|\leq\frac{1}{\sqrt{N}}\|\Lambda\|\|\hSig_u^{-1}-\Sigu^{-1}\|\|u_t\|=O_p(\sqrt{N}m_N\omega_T^{1-q}),
 $$
 which is $O_p(1+\sqrt{N(\log N)/T})$ even if $m_N$ is bounded and $q=0$. Hence this leads to a crude bound that does not converge.  The problem is present even if    $\Sigma_u^{-1}$ is estimated with the optimal rate of convergence.  

We realize that such a technical problem is common for statistical inferences that involve estimating a high-dimensional covariance.     In fact, most of the existing  approaches in the literature   only produce ``absolute convergence" $\|\widehat\Sigma_u^{-1}-\Sigma_u^{-1}\|$. For statistical inference purposes like the primary interest of this paper, however,  the absolute convergence is  not sufficient when $N/T\rightarrow\infty$.

 We propose a new technical strategy to solve this problem, by directly investigating the  ``weighted convergence" of the  weighted error:
\begin{equation}\label{eq4.6}
\|\frac{1}{\sqrt{N}}\Lambda'(\hSig_u^{-1}-\Sigma_u^{-1})u_t\|.
\end{equation}
Intuitively, the weights $\Lambda'$ and $u_t$ ``average down" the estimation errors, and improve the rate of convergence.
Formal analysis  requires us to re-investigate the asymptotic behavior of the thresholded covariance estimator.  We require the following technical assumption.


Let  $\Lambda'\Sigma_u^{-1}=(\xi_1,...,\xi_N)$. Assuming $\|\Sigma_u^{-1}\|_1=O(1)$, we then have  $\max_{j\leq N}\|\xi_j\|<C$ for some $C>0.$ In addition, let $e_t=\Sigma_u^{-1}u_t$, then $e_t$ has mean zero and covariance $\Sigma_u^{-1}$. 

\begin{assum}\label{a43} For each $t\leq T$ and $k\leq r$,\\
(i) $\frac{1}{T\sqrt{N}}\sum_{i=1}^N\sum_{s=1}^T(u_{is}^2-Eu_{is}^2)\xi_ie_{it}=o_p(1)$\\
(ii) $\frac{1}{NT\sqrt{N}}\sum_{i=1}^N\sum_{s=1}^T\sum_{j=1}^N (u_{js}u_{is}-Eu_{js}u_{is})\lambda_j\lambda_i'e_{it}\xi_{ik}=o_p(1)$, \\
(iii)$\frac{1}{T\sqrt{N}}\sum_{i\neq j, (i,j)\in S_U}\sum_{s=1}^T(u_{is}u_{js}-Eu_{is}u_{js})\xi_ie_{jt}=o_p(1)$,\\
(iv) $\frac{1}{NT\sqrt{N}}\sum_{i\neq j, (i,j)\in S_U}\sum_{v=1}^N\sum_{s=1}^T(u_{is}u_{vs}-Eu_{is}u_{vs})      \xi_{ik}e_{jt}\lambda_{v}\lambda_{j}'=o_p(1)$.

\end{assum}
The above conditions are new in the literature and essential to establish the  weighed convergence.   The intuition of these conditions  is that, the weighted average of the standardized sum $\frac{1}{\sqrt{T}}\sum_{t=1}^T(u_{it}u_{jt}-Eu_{it}u_{jt})$  is   $o_p(1)$ once averaged across $i$ and $j$.    The extra term $\frac{1}{N}$ appeared in $\frac{1}{NT\sqrt{N}}$ of  Conditions (ii) and (iv) is a scaling factor because under the sparsity condition, the number of summands of $\sum_{i=1}^N$ and  $\sum_{i\neq j, (i,j)\in S_U}$   is at most $O(N)$ (e.g., in block diagonal and banded matrices).  

We verify the key assumption \ref{a43}  in  the following lemma, when  $\{u_t\}_{t\leq T}$ is serially independent. We require  $N=o(T^2)$ but still allow $N/T\rightarrow\infty$.

\begin{lem}\label{l4.1}
Suppose $\{u_{it}\}_{t\leq T}$  is  independent across $t$ (but can still be correlated across $i$), and the sparse condition Assumption \ref{a41} holds.  Then when $N=o(T^2)$, Assumption \ref{a43} is satisfied.
\end{lem}

We have the following \textit{weighted consistency} for the estimated weight matrix, which as we have explained, cannot be implied directly by the absolute convergence $\|\widehat\Sigma_u^{-1}-\Sigma_u^{-1}\|$ even when $\Sigma_u$ is diagonal. As one of the main contributions of this paper, result of this type is potentially widely useful for high-dimensional inferences when large covariance estimation is involved.

\begin{prop}\label{l41}  Suppose   $\sqrt{N}m_N^2\omega_T^{2-2q}=o(1)$, and Assumptions  \ref{a21}-  \ref{a28},  
  \ref{a41}, \ref{a43} hold. For $q$, $m_N$ and $\omega_T$  defined in (\ref{e40}) and (\ref{e42}), and for each $t\leq T$,     we have
  $$
\|\frac{1}{\sqrt{N}}\Lambda'(\hSig_u^{-1}-\Sigu^{-1})u_t\|=o_p(1).
$$
Therefore  Assumption \ref{a26} is satisfied for $W=\Sigma_u^{-1}$.
\end{prop}

\begin{remark}
Consider a strictly sparse case where \\$m_N=\max_{i\leq N}\sum_{j=1}^NI(\Sigma_{u,ij}\neq0)=O(1)$. The condition in the theorem  $\sqrt{N}m_N^2\omega_T^{2-2q}=o(1)$ then holds as long as $\sqrt{N}\log N=o(T)$. As  always the case,  requiring $N=o(T^2)$ is needed for the asymptotic normality of  $\widehat f_t.$

\end{remark}


 \subsection{Efficient estimation}
We use $W_T=\hSig_u^{-1}$ as the feasible  weight matrix.  Let  the columns of the $T\times r$ matrix $\hF^e/\sqrt{T}=(\hf_1^e,...,\hf_T^e)'/\sqrt{T}$ be  the eigenvectors corresponding to the largest $r$ eigenvalues of $Y'\hSig_u^{-1}Y$, and $\hLam^e=T^{-1}Y\hF^e=(\hlam_1^e,...,\hlam_N^e)'.$   Here the superscript $e$ denotes ``efficient" WPC.  


We denote $\Sigma_{\Lambda,e}$ as the limit of $\Lambda'\Sigma_u^{-1}\Lambda/N$. Let $V_e$ be an $r\times r$ diagonal matrix with elements as the largest $r$ eigenvalues of $\Sigma_{\Lambda, e}^{1/2}\cov(f_t)\Sigma_{\Lambda, e}^{1/2}$, and  $\Gamma_e$ be the corresponding eigenvector matrix such that $\Gamma_e'\Gamma_e=I_r$.   In addition,  let $Q_e=V_e^{1/2}\Gamma_e'\Sigma_{\Lambda, e}^{-1/2}$.   We have the following limiting distributions for the estimated factors and loadings.
\begin{thm} \label{t41} Under the assumptions of Proposition \ref{l41}, for each $t\leq T$ and $j\leq N$,
$$
\sqrt{T}(\hlam_j^e-H_e^{'-1}\lambda_j)\rightarrow^d \mathcal{N}(0, Q_e^{'-1}\Phi_jQ_e^{-1}).
$$
$$
\sqrt{N}(\hf_t^e-H_ef_t)\rightarrow^d \mathcal{N}(0, V_e^{-1}).
$$where $\Phi_j$ is as defined in (\ref{eq2.9}). In addition, for the estimated common component,
 $$
\frac{\hlam_i^{e'}\hf_t^e-\lambda_i'f_t}{(\lambda_i'\Xi_e\lambda_i/N+f_t'\Omega_if_t/T)^{1/2}}\rightarrow^d\mathcal{N}(0,1).
$$
where $\Xi_e=(\Lambda'\Sigma_u^{-1}\Lambda/N)^{-1}
$ and $\Omega_i$ is defined as in Theorem \ref{t33}.
\end{thm}

For completeness, the following result gives  the uniform rate of convergence.
\begin{thm} \label{c41} Suppose $N^{1/(2-2q)}\log N=o(T)$ and $T=o(N^2)$. Under the assumptions of Theorem \ref{l41}, there is an $r\times r$ matrix $H_e$ such that
\begin{eqnarray*}
\max_{t\leq T}\|\hf_t^e-H_ef_t\|&=&O_p\left(\frac{T^{1/\delta}}{\sqrt{N}}+(\log T)^{\alpha}m_N\omega_T^{1-q}\right),\\
\max_{j\leq N}\|\hlam_j^e-H_e^{'-1}\lambda_j\|&=&O_p\left(m_N\omega_T^{1-q}\right).
\end{eqnarray*}
\end{thm}

\begin{remark}
Typically in the  strictly sparse case  $m_N=O(1)$ and $q=0$. When $N/T\rightarrow\infty$, the above rates become:
$$
\max_{t\leq T}\|\hf_t^e-H_ef_t\|=O_p\left(\frac{T^{1/\delta}}{\sqrt{N}}+\frac{(\log T)^{\alpha}\sqrt{\log N}}{\sqrt{T}}\right),$$
$$ 
\max_{j\leq N}\|\hlam_j^e-H_e^{'-1}\lambda_j\|=O_p\left(\sqrt{\frac{\log N}{T}}\right).
$$

\end{remark}



\subsection{Optimal weight matrix}\label{sow}

Regular PC, heteroskedastic WPC and   the efficient WPC minimize  different objective functions, depending  on the  choices of the weight matrix. Thus the estimated $\hF/\sqrt{T}$ are the eigenvectors  from three different matrices. Table \ref{table1} summarizes the main differences of the estimators.

\begin{table}[htdp]
\caption{Three interesting choices of $W$ }
\begin{center}
\begin{tabular}{c|c|c|c}
\hline
&Objective function& Eigenvectors of & $W$\\
\hline
regular PC & $\sum_{t=1}^T(Y_t-\Lambda f_t)'(Y_t-\Lambda f_t)$ &  $Y'Y  $ & $I_r$\\
&&\\
HWPC & $\sum_{t=1}^T(Y_t-\Lambda f_t)'\diag(\hSig_u)^{-1}(Y_t-\Lambda f_t)$ & $Y'\diag(\hSig_u)^{-1}Y   $ &$\diag(\Sigma_u)^{-1}$\\
&&\\
EWPC & $\sum_{t=1}^T(Y_t-\Lambda f_t)'\hSig_u^{-1}(Y_t-\Lambda f_t)$&  $Y'\hSig_u^{-1}Y $& $\Sigu^{-1}$\\
\hline
\end{tabular}
\label{table1}

\end{center}
\small\it  The estimated $\hF/\sqrt{T}$  is the eigenvectors of the largest $r$ eigenvalues of $Y'W_TY$, and $\hLam=T^{-1}Y\hF$. 
HWPC represents the heteroskedastic WPC; EWPC  represents the efficient WPC.
\end{table}


A natural question arises: is the consistent estimator for $W=\Sigu^{-1}$ indeed the optimal choice over  a broad class of positive definite weight matrices? One can answer this question via looking at the asymptotic variance of the estimators, as choosing the optimal weight for GMM (Hansen 1982). However, because   WPC estimators are  estimating rotated factors and loadings,   the rotation   depends on the choice of $W.$ 
But regardless of the choice $W$, the common component $\lambda_i'f_t$ is always directly estimated.  
The following result demonstrates that $W_T=\widehat\Sigma_u^{-1}$   yields the minimum asymptotic variance of $\hlam_i'\hf_t$ among WPC estimators.

\begin{thm}\label{t43} Let $(\lambda_i'\Xi_e\lambda_i/N+f_t'\Omega_if_t/T)$ denote the asymptotic variance of $\hlam^{e'}_i\hf_t^e$ based on $\hSig_u^{-1}$ as in Theorem \ref{t41}. For any positive definite matrix $W$, let $(\lambda_i'\Xi_W\lambda_i/N+f_t'\Omega_if_t/T)$ denote the asymptotic variance of $\hlam'_i\hf_t$ as in Theorem \ref{t33} based on $W$.  Then for each  $i\leq N$ and $t\leq T$,
$$
\lambda_i'\Xi_e\lambda_i/N+f_t'\Omega_if_t/T\leq \lambda_i'\Xi_W\lambda_i/N+f_t'\Omega_if_t/T.
$$
In fact, for all large $N$,
$
\Xi_W-\Xi_e
$
is semi-positive definite for each positive definite matrix $W$.  
\end{thm}




\subsection{Estimating asymptotic covariances}\label{s45}

We derive  consistent estimators for the asymptotic  variances  that appeared in Theorem \ref{t41}. Hence the derived optimal limiting distributions can be used for statistical inferences. These estimators account for the serial and cross-sectional correlations of the data in both $i$ and $t$.

The factor estimator has an asymptotic expansion:
$$
\sqrt{N}(\hf_t^e-H_ef_t)=\hV^{-1}\frac{\hF^{e'}F}{T}\frac{\Lambda'\Sigma_u^{-1}u_t}{\sqrt{N}}+o_p(1)
$$
where  $\hV$ is the $r\times r$ diagonal matrix of the first $r$ largest eigenvalues of $\frac{1}{TN}Y\hSig_u^{-1} Y'.$    Theorem \ref{t41} shows that the asymptotic variance is $V_e^{-1}$. Hence,
\begin{equation}\label{e45} 
\hV^{-1}\frac{\hF^{e'}F}{T}\frac{\Lambda'\Sigma_u^{-1}\Lambda}{N}\frac{F'\hF^e}{T}\hV^{-1} \rightarrow^p  V_e^{-1}
\end{equation}
The left hand side involves the product $F\Lambda'$, which can be estimated by $\hF^e\hLam^{e'}.$  A
consistent estimator of $V_e^{-1}$ is then given by (note that $\frac{1}{T}\hF^{e'}\hF^e=I_r$)
$$
\hV_e^{-1}=\hV^{-1}\frac{\hF^{e'}\hF^e}{T}\frac{\hLam^{e'}\hSig_u^{-1}\hLam^e}{N}\frac{\hF^{e'}\hF^e}{T}\hV^{-1} =\frac{1}{N}\hV^{-1}\hLam^{e'}\hSig_u^{-1}\hLam^e\hV^{-1}.
$$

The loading estimator has an asymptotic expansion:
$$
\sqrt{T}(\hlam_j-H_e^{'-1}\lambda_j)=\frac{1}{\sqrt{T}}\sum_{t=1}^TH_ef_tu_{jt}+o_p(1).
$$
Here  $H_ef_tu_{jt}$ can be   estimated by $\hf_t^e\hu_{jt}$, where $\hu_{jt}$ is a WPC estimator of the error term (e.g.,  $\hu_{jt}=y_{it}-\hlam_j^{e'}\hf_t^e$). We apply the HAC (heteroskedasticity
and autocorrelation consistent) estimator of Newey and West (1987) to estimate $Q_e^{'-1}\Phi_jQ_e^{-1}$, the asymptotic variance of $\sqrt{T}(\hlam_j-H_e^{'-1}\lambda_j)$, based on the sequence $\{\hf_t^e\hu_{jt}\}$:
$$
\widehat\Psi_j=\frac{1}{T}\sum_{t=1}^T\hu_{jt}^2\hf_t^e\hf_t^{e'}+\sum_{l=1}^K(1-\frac{l}{K+1})\frac{1}{T}\sum_{t=l+1}^T\hu_{jt}\hu_{j,t-l}(\hf_t^e\hf_{t-l}^{e'}+\hf_{t-l}^e\hf_{t}^{e'}),
$$
where $K=K_{T, N}\rightarrow\infty$ is an increasing sequence such that \\$K=o(\min\{T^{1/4}, N^{1/4}\}).$ The advantages of using the HAC estimator are:   it accounts for the serial correlations of $\{f_tu_t'\}_{t\geq 1}$, and it  also guarantees the positive semi-definiteness for any given finite sample as shown by Newey and West (1987).

The asymptotic variance of the common component in Theorem \ref{t41} consists of $\lambda_i'\Xi_e\lambda_i$ and $f_t'\Omega_if_t$, where $\Xi_e=(\frac{1}{N}\Lambda'\Sigma_u^{-1}\Lambda)^{-1}$ and \\$\Omega_i=\cov(f_t)^{-1}\Phi_i\cov(f_t)^{-1}$. We respectively estimate them by $$
\widehat\Theta_{1i}=\frac{1}{N}\hlam_i^{e'}\hV^{-1}\hLam^{e'}\hSig_u^{-1}\hLam^e\hV^{-1}\hlam_i^e,
\quad 
\widehat\Theta_{2, it}=\hf_t^{e'}\widehat\Psi_i\hf_t^e.
$$
\begin{thm}\label{t44} Under the assumptions of Theorem  \ref{c41}, as $T$, $N\rightarrow\infty$,  and $K=K_{T, N}=o(\min\{T^{1/4}, N^{1/4}\})$,
\begin{eqnarray*}
&&\hV_e^{-1}\rightarrow^p  V_e^{-1},  \qquad \widehat\Psi_j\rightarrow^p Q_e^{'-1}\Phi_jQ_e^{-1},\cr
&&\widehat\Theta_{1i}\rightarrow^p\lambda_i'\Xi_e\lambda_i, \qquad 
\widehat\Theta_{2, it}\rightarrow^pf_t'\Omega_i f_t.
\end{eqnarray*}
\end{thm}

These covariance estimators can be easily computed. 

\section{WPC for Panel data Models with Interactive Effects}\label{secIn}
The factor model we have considered so far is closely related to the following panel data model:
\begin{equation}\label{e51}
y_{it}=X_{it}'\beta+\varepsilon_{it},\quad \varepsilon_{it}=\lambda_i'f_t+u_{it},\quad i\leq N, t\leq T
\end{equation}
 The regression noise   has a factor structure with unknown $\lambda_i$ and $f_t$, and $u_{it}$ still represents  the idiosyncratic error component.  It is assumed that $u_{it}$ is independent of $(X_{it}, f_t)$.   In the model, the only observables are $(y_{it}, X_{it})$. The goal is to estimate $\beta$, the structural parameter of the model. 

Substituting the second equation to the first one in (\ref{e51}), we obtain
\begin{equation}\label{e52}
y_{it}=X_{it}'\beta+\lambda_i'f_t+u_{it}.
\end{equation}
If we treat $\lambda_i$ as the ``individual effect" and $f_t$ as the ``time effect", then the factor structure $\lambda_i'f_t$ represents the interaction between the individual and time effects, so called ``interactive effect".  This model was previously studied by, e.g., Ahn et al. (2001),  Pesaran (2006), Bai (2009), Moon and Weidner (2010). 

The difficulty of estimating $\beta$ is that, in many applied problems  the regressor $X_{it}$ is correlated with  the  time effect (common factor) $f_t$, individual effect $\lambda_i$, or both.  As a result, $X_{it}$ and $\varepsilon_{it}$ are also correlated, so regressing $y_{it}$ directly on $X_{it}$ cannot produce a consistent estimator for $\beta$.  In addition, existing methods ignore the  heteroskedasticity and correlation in $\{u_{it}\}_{i\leq N}$. Hence efficiency is lost, for instance, when   $\Sigma_u$ is non-diagonal or its diagonal entries vary over a large range. We shall illustrate the consequence of  efficiency loss using a real data application in Section \ref{law}.



\subsection{WPC estimation of $\beta$}  Let $X_t=(X_{1t},...,X_{Nt})'$,  $(N\times d)$. 
We estimate $\beta$ via 
\begin{equation}\label{eq6.3add}
\min_{\beta, f_t, \Lambda}\sum_{t=1}^T(Y_t-\Lambda f_t-X_t\beta)'W(Y_t-\Lambda f_t-X_t\beta),
\end{equation}
for some positive definite $N\times N$ weight matrix. Similar to the  generalized least squares estimator (GLS) for linear regressions,   we choose the weight matrix to be
$$
W=\Sigma_u^{-1}.
$$
This choice produces   similar estimators as the efficient WPC. The estimator is  feasible once we consistently estimate $\Sigma_u^{-1}$, which can be done under the assumption that $\Sigma_u$ is sparse. Suppose $\widetilde{\Sigma}_u^{-1}$ is a consistent covariance estimator. The feasible WPC estimates   $\beta$ by:
\begin{equation}
\hat\beta=\arg\min_{\beta}\min_{f_t, \Lambda}\sum_{t=1}^T(Y_t-\Lambda f_t-X_t\beta)'\tSig_u^{-1}(Y_t-\Lambda f_t-X_t\beta),
\end{equation}
where the minimization is subjected to the constraint $\frac{1}{T}\sum_{t=1}^Tf_tf_t'/T=I_r$ and $\Lambda'\tSig_u^{-1}\Lambda$ being diagonal. The estimated $\beta$ for each given $(\Lambda, f_t)$ is simply
$$
\beta(\Lambda, f_t)=(\sum_{t=1}^TX_t'\tSig_u^{-1}X_t)^{-1}\sum_{t=1}^TX_t'\tSig_u^{-1}(Y_t-\Lambda f_t).
$$
On the other hand, given $\beta$, the variable $Y_t-X_t\beta$ has a factor structure. Hence the estimated $(\Lambda, f_t)$ are the WPC estimators: let $X(\hat\beta)$ be an $N\times T$ matrix $X(\hat\beta)=(X_1\hat\beta,...,X_T\hat\beta).$ The columns of the $T\times r$ matrix $\tF/\sqrt{T}=(\tf_1,...,\tf_T)'/\sqrt{T}$ are  the eigenvectors corresponding to the largest $r$ eigenvalues of $(Y-X(\hat\beta))'\tSig_u^{-1}(Y-X(\hat\beta))$, and $\widetilde{\Lambda}=T^{-1}(Y-X(\hat\beta))\tF.$ Therefore, given $(\Lambda, f_t)$, we can estimate $\beta$, and given $\beta$, we can estimate $(\Lambda, f_t)$.  So $\hat\beta$ can be simply obtained by iterations, with an initial value $\hat\beta_0$. This iteration scheme only requires two matrix inverses: $\tSig_u^{-1}$ and  $(\sum_{t=1}^TX_t'\tSig_u^{-1}X_t)^{-1}$, which do not update during iterations.  
Based on our experience of numerical studies, the  iterations converge   fast.

Similar to Fan et al. (2013), the covariance estimator can be constructed based on thresholding. Let $\hat\beta_0$ be a   ``regular PC  estimator" that takes $W=I_N$ in (\ref{eq6.3add}), 
  which is  known to be $\sqrt{NT}$-consistent (e.g., Bai 2009, Moon and Weidner 2010).  Apply the singular value decomposition to
$$
\frac{1}{T}\sum_{t=1}^T(Y_t-X_t\hat\beta_0)(Y_t-X_t\hat\beta_0)'=\sum_{i=1}^N\nu_ig_ig_i',
$$
where $(\nu_j,g_j)_{j=1}^N$ are the  eigenvalues-eigenvectors of $\frac{1}{T}\sum_{t=1}^T(Y_t-X_t\hat\beta_0)(Y_t-X_t\hat\beta_0)'$  in a decreasing order such that $\nu_1\geq \nu_2\geq...\geq\nu_N.$   Then $\widetilde{\Sigma}_u=(\widetilde{\Sigma}_{u,ij})_{N\times N}$,
$$
\tSig_{u,ij}=\begin{cases}
\widetilde R_{ii}, & i=j\\
s_{ij}(\widetilde R_{ij}), & i\neq j
\end{cases}, \quad \widetilde R=(\widetilde R_{ij})_{N\times N}=\sum_{i=r+1}^N\nu_ig_ig_i',
$$
where $s_{ij}(\cdot)$ is the same thresholding function  as defined in Section \ref{swc} with the same threshold $\tau_{ij}$.


\subsection{Assumptions for asymptotic analysis}\label{s5.2}
     Rearrange the design matrix
$$
Z=(X_{11},..., X_{1T}, X_{21},...,X_{2T},...,X_{N1},...,X_{NT})',\quad NT\times d.
$$
For any $T\times r$ matrix $F$, let $M_F=I_T-F(F'F)^{-1}F'/T$. The following matrices play an important role in the identification and asymptotic analysis:
\begin{eqnarray}\label{eq5.6}
A_F&=&\left[\Sigma_u^{-1}-	\Sigma_u^{-1}\Lambda\left(\Lambda'\Sigma_u^{-1}\Lambda\right)^{-1}\Lambda'\Sigma_u^{-1}		\right]\otimes M_F, 	\cr
V(F)&=&\frac{1}{NT}Z' A_FZ,
\end{eqnarray}
where $(\Lambda,\Sigma_u^{-1})$ in the above represent the true loading matrix and inverse error covariance in the data generating process, and $\otimes$ denotes the Kronecker product. Our first condition assumes that $V(F)$ is positive definite in the limit uniformly over a class of  $F$.
\begin{assum}\label{a51}
With probability approaching one,
$$
\inf_{F: F'F/T=I_r} \lambda_{\min}(V(F))>0.
$$
\end{assum}
If we write $B_F=\left[\Sigma_u^{-1/2}-	\Sigma_u^{-1}\Lambda\left(\Lambda'\Sigma_u^{-1}\Lambda\right)^{-1}\Lambda'\Sigma_u^{-1/2}		\right]\otimes M_F, 	$ then $A_F=B_FB_F'.$ So $V(F)$ is at least semi-positive definite. Also,  summing over $NT$ rows of $Z$ should lead to a strictly positive definite matrix $V(F)$. 
As a sufficient condition, if $X_{it}$  depends on the factors and loadings through:
$$
X_{it}=\tau_i+\theta_t+\sum_{k=1}^ra_k\lambda_{ik}+\sum_{k=1}^rb_kf_{kt}+\sum_{k=1}^rc_k\lambda_{ik}f_{kt}+\eta_{it}
$$
where $a_k, b_k, c_k$ are constants (can be zero) and $\eta_{it}$ is i.i.d. over both $i$ and $t$,  then Assumption \ref{a51} is satisfied (see Bai 2009). 

 Let $U=(u_{11},...,u_{1T},u_{21},..., u_{2T},..., u_{N1},...,u_{NT})'$, and   $F_0$ be  the $T\times r$ matrix of true factors.
\begin{assum}\label{ass5.3}  There is a $\dim(\beta)\times \dim(\beta)$  positive definite matrix $\Gamma$ such that
\begin{eqnarray*}
V(F_0)\rightarrow^p\Gamma,\quad\frac{1}{\sqrt{NT}}Z'A_{F_0}U\rightarrow^d\mathcal{N}(0,\Gamma).
\end{eqnarray*}
\end{assum}
This assumption is required for the asymptotic normality of $\hat\beta$, because it can be shown that,  
$$
\sqrt{NT}(\hat\beta-\beta) =V(F_0)^{-1}\frac{1}{\sqrt{NT}}Z'A_{F_0}U+o_p(1).
$$
Hence the asymptotic normality depends on that of $\frac{1}{\sqrt{NT}}Z'A_{F_0}U$.
Assumption \ref{ass5.3} is not stringent because  if we write $B_{F_0}'U=(\widetilde{u}_{11},...,\widetilde{u}_{1T}, \widetilde{u}_{21},...,\widetilde{u}_{NT})' $, and $Z'B_{F_0}=(\widetilde{Z}_{11},...,\widetilde{Z}_{1T}, \widetilde{Z}_{21},...,\widetilde{Z}_{NT})$, then\\
$
\frac{1}{\sqrt{NT}}Z'A_{F_0}U=\frac{1}{\sqrt{NT}}\sum_{t=1}^T\sum_{i=1}^N\widetilde{Z}_{it}\widetilde{u}_{it}
$ is a standardized summation. We can further write
$$
\sqrt{NT}(\hat\beta-\beta)=\left(\frac{1}{NT}\sum_{t=1}^T\sum_{i=1}^N\widetilde{Z}_{it}\widetilde{Z}_{it}'\right)^{-1}\frac{1}{\sqrt{NT}}\sum_{t=1}^T\sum_{i=1}^N\widetilde{Z}_{it}\widetilde{u}_{it}+o_p(1).$$
Hence the second statement of Assumption \ref{ass5.3} is  a central limit theorem for $\frac{1}{\sqrt{NT}}\sum_{t=1}^T\sum_{i=1}^N\widetilde{Z}_{it}\widetilde{u}_{it}$  on both cross-sectional and time domains. 
 In addition, in the absence of serial correlation, the conditional covariance of $\frac{1}{\sqrt{NT}}Z'A_{F_0}U$ given $Z$ and $F_0$ equals $\frac{1}{NT}Z'A_{F_0}(\Sigma_u\otimes I_T)A_{F_0}Z=V(F_0)$. This implies that the asymptotic variance of $\sqrt{NT}(\hat\beta-\beta_0)$ is simply $\Gamma^{-1}$.

 
 \subsection{Weighted convergence for estimating the weight matrix}\label{eew}
 
 The issue described in Section 2 arises in establishing
  \begin{equation}\label{eq6.8}
\frac{1}{\sqrt{NT}}Z'[(\tSig_u^{-1}-\Sigma_u^{-1})\otimes I_T]U=o_p(1),
 \end{equation}
 which is the effect of estimating the large covariance $\Sigma_u^{-1}$.
 In fact, the first order condition of $\hat\beta$ leads to 
 $$
 \sqrt{NT}(\hat\beta-\beta)=V(F_0)^{-1}\frac{1}{\sqrt{NT}}Z'\widehat AU+o_p(1),
 $$
 where $\widehat A$ is as $A_{F_0}$ with $\Sigma_u^{-1}$ replaced with $\tSig_u^{-1}$ and $F_0$ replaced with $\widetilde F$. Hence we need 
 \begin{equation}\label{eq5.7}
 \frac{1}{\sqrt{NT}}Z'(\widehat A-A_{F_0})U=o_p(1).
 \end{equation}
This requires   the weighted convergence (\ref{eq6.8}).  However, when $N/T\rightarrow\infty$,  achieving (\ref{eq6.8}) is technically  difficult.  Similar to the case described in the approximate factor model,  the absolute convergence of $\|\tSig_u^{-1}-\Sigma_u^{-1}\|$  is not suitable for inferences.



We consider the Gaussian case for simplicty, and the problem is still highly technically involved. Non-Gaussian case will be even more challenging, and we shall leave it for future research. 
\begin{assum}\label{a53add}

(i)  $u_{t}$ is distributed as $\mathcal{N}(0,\Sigma_u)$.  \\
(ii) $\{u_{t}\}_{t\geq 1}$ is independent of $\{f_t, X_t\}_{t\geq 1}$, and $\{u_t, f_t, X_t\}$ are serially independent across $t$.
  
\end{assum}

 It is possible to relax Condition (ii) to allow for serial correlations, but $\hat\beta$ will be asymptotically biased.

\subsection{Limiting distribution}
We  require the same conditions on the data generating process for  the factors,  loadings and the sparsity of $\Sigma_u$ as in Sections 2 and \ref{sepc}.

\begin{prop}\label{prop2}
Under Assumptions \ref{a21}- \ref{a27},  \ref{a41},  \ref{a51}-\ref{a53add}, as $N/T\rightarrow\infty$, and $m_N=o(T^2)$,  we have the weighted convergence:
$$ \frac{1}{\sqrt{NT}}Z'(\widehat A-A_{F_0})U=o_p(1).$$
\end{prop}

We have the following limiting distribution.

\begin{thm} \label{t51} Under  the assumptions of Proposition \ref{prop2}, the asymptotic limiting distribution of $\hat\beta$ is the same when either   $W=\Sigma_u^{-1}$ or  the feasible weight $W_T=\tSig_u^{-1}$  is used as the weight matrix,  and is given by
$$
\sqrt{NT}(\hat\beta-\beta)\rightarrow^d\mathcal{N}(0,  \Gamma^{-1}).
$$
\end{thm}

The asymptotic variance $\Gamma^{-1}$ is the limit of $V(F_0)^{-1}$.   Note that under the same set of conditions, the regular PC method of Bai (2009) and Moon and Weidner (2010) gives an  asymptotic conditional covariance (given $Z, F_0$) of the sandwich-formula: $$V_2\equiv (\frac{1}{NT}Z'GZ)^{-1}\frac{1}{NT}Z'G(\Sigma_u\otimes I_T)GZ(\frac{1}{NT}Z'GZ)^{-1},$$ where $G$ is defined as $A_{F_0}$ with $\Sigma_u^{-1}$ replaced with $I_N.$ It is not hard to show that $V_2-V(F_0)^{-1}$ is semi-positive definite. So relative efficiency is gained when  WPC is used. In fact, the choice $W=\widetilde\Sigma_u^{-1}$ is also the optimal weight matrix for WPC in this case.


To estimate the asymptotic variance of $\hat\beta$, let $\widetilde A$   equal  $A_{F}$ with $F$, $\Lambda$ and $\Sigma_u^{-1}$ replaced with $\widetilde F$, $\widetilde\Lambda$ and $\widetilde\Sigma_u^{-1}$. Define $\widetilde\Gamma=\frac{1}{NT}Z'\widetilde AZ$.  The following result enables us to construct confidence intervals and conduct hypothesis tests for $\beta$ under large samples.

\begin{thm}
\label{t52} Under the assumptions of Theorem  \ref{t51}, 
$$
\widetilde\Gamma^{-1}\rightarrow^p \Gamma^{-1}.
$$
\end{thm}

The methods of Section \ref{sepc} also  carry  over to derive the limiting distributions of the estimated interactive effects $\lambda_i'f_t$. The procedure and corresponding results are very similar given  the  $\sqrt{NT}$-consistency of $\hat\beta$. Hence we omit  repeated discussions. 


 \subsection{Estimation with unknown number of factors} 
 
 For simplicity of presentations, we have assumed the number of factors $r$ to be known. As was shown by many authors,  estimation results are often robust to over-estimating $r$. For instance, Moon and Weidner  (2011) have shown that for inference on the regression coefficients one does not need to estimate $r$ consistently, as long as the ``working number" is not less than the true value. On the other hand, we can also start with a consistent estimator $\hat r$ using a similar method of Bai and Ng (2002) and Bai (2009). 
 
 Specifically, suppose there is a known upper bound $\bar r$ of the number of factors. For each $k\leq \bar r$, define 
 $$
 \widehat\sigma^2(k)=\min_{\beta,\Lambda_k, f_{t,k}}\frac{1}{NT}\sum_{t=1}^T(Y_t-\Lambda_k'f_{t,k}-X_{t}\beta)'(Y_t-\Lambda_k'f_{t,k}-X_{t}\beta) $$
 where each row of $\Lambda_k$ is  a $k$-dimensional loading vector, and $f_{t,k}$ is also $k$-dimensional. The above minimization is subject to the constraint that $\frac{1}{T}\sum_{t=1}^Tf_{t,k}f_{t,k}'=I_k$ and $\Lambda_k'\Lambda_k$ is diagonal. The iterative algorithm based on principal components  can calculate the above minimization fast. Under our conditions, Bai (2009) showed that $r$ can be consistently estimated by either minimizing CP$(k)$ or IC$(k)$, where
 $$
 \text{CP}(k)=\widehat\sigma^2(k)+\widehat\sigma^2(\bar k)[k(N+T)-k^2]\frac{\log (NT)}{NT},
 $$
 and
 $$
  \text{IC}(k)=\log \widehat\sigma^2(k)+[k(N+T)-k^2]\frac{\log (NT)}{NT}.
 $$
 We then can apply the estimator $\hat r$ to construct the WPC estimator, and achieve the same limiting distributions. Estimation procedure and its theoretical properties can be proved to be the same as before, so details are not presented to avoid repetition.

\section{Simulated Experiments}\label{ssim} 
We conduct numerical experiments to compare  the proposed WPC with the popular methods in the literature\footnote{
We have written a Matlab code to   implement   the proposed  WPC for any user-specified weight matrix as well as the optimal WPC for both the factor model and panel data model with interactive effects, available upon request.  }. The idiosyncratic error terms are  generated   as follows:  let  $\{\epsilon_{it}\}_{i\leq N, t\leq T}$ be  i.i.d.  $\mathcal{N}(0,1)$ in both $t,i$. Let
$$u_{1t}=\epsilon_{1t}, \hspace{1em}u_{2t}=\epsilon_{2t}+a_1\epsilon_{1t},  \hspace{1em}u_{3t}=\epsilon_{3t}+a_2\epsilon_{2t}+b_1\epsilon_{1t},$$
$$
u_{i+1,t}=\epsilon_{i+1,t}+a_i\epsilon_{it}+b_{i-1}\epsilon_{i-1,t}+c_{i-2}\epsilon_{i-2,t},
$$
where $\{a_i, b_i, c_i\}_{i=1}^N$ are  i.i.d. $\mathcal{N}(0, 1)$.  Then $\Sigma_{u}$ is a banded matrix, possessing both cross-sectional correlation and heteroskedasticity. Let the two factors $\{f_{1t}, f_{2t}\}$ be i.i.d. $\mathcal{N}(0,1)$, and $\{\lambda_{i,1}, \lambda_{i,2}\}_{i\leq N}$ be uniform on $[0,1]$.  We estimate the optimal weight matrix by soft-thresholding the ``correlation matrix" of $R$ as suggested by Fan et al. (2013). 

 \textbf{Design 1}
 
 Consider the pure factor model $y_{it}=\lambda_{i1}f_{1,t}+\lambda_{i,2}f_{2t}+u_{it}$, where we estimate  the factor loadings $\{\lambda_{i,1}, \lambda_{i,2}\}_{i\leq N}$ and factors $\{f_{1t}, f_{2t}\}$. For each estimator, the smallest canonical correlation (the larger the better) between the estimators and   parameters are calculated,  as  an assessment of the estimation  accuracy.   The simulation is  replicated for one hundred times, and the average canonical correlations for several competing methods are reported in Table \ref{table2}.  
 The mean squared error of the estimated common components are also compared.

  \begin{table}[htdp]
\caption{Canonical correlations  for simulation study}
\begin{center}
\begin{tabular}{cc|ccc|ccc|ccc}
 
\hline
&& \multicolumn{3}{c|}{ Loadings}  &\multicolumn{3}{c|}{ Factors} &\multicolumn{3}{c}{$(\frac{1}{NT}\sum_{i,t}(\widehat\lambda_i'\widehat f_t-\lambda_i'f_t)^2)^{1/2}$}  \\
$T$ & $N$ &PC& HWPC & EWPC &PC&HWPC &EWPC&PC&HWPC &EWPC\\
 &   &\multicolumn{3}{c|}{ (the larger the better)}  &\multicolumn{3}{c|}{ (the larger the better)} &\multicolumn{3}{c}{ (the smaller the better)} \\
\hline
& & & & & &  &&\\
50 & 75 & 0.346&0.429 &  0.487&  0.403   &  0.508&   0.566& 0.621&0.583&0.545\\
50 & 100 &   0.411&   0.508 &   0.553&  0.476    & 0.602     &  0.666&0.546&0.524&0.498\\
50 & 150 &  0.522&  0.561 &  0.602 &  0.611&  0.679&     0.746&0.467&0.444&0.427\\
& & & & & & &&\\
100 & 80 & 0.433& 0.545 &   0.631&  0.427   & 0.551&  0.652&0.570&0.540&0.496\\
100 & 150 &  0.613&   0.761&   0.807 &  0.661    &  0.835&   0.902&0.385& 0.346&  0.307\\
100 & 200 &   0.751&  0.797 &    0.822   & 0.827& 0.882& 0.924&0.333&0.312&0.284\\
& & & & & &  &&\\
150 & 100 &   0.380&  0.558&  0.738 &   0.371   & 0.557  &  0.749&0.443& 0.394& 0.334\\
150 & 200 &   0.836&  0.865&  0.885 & 0.853    &  0.897&  0.942&0.313&0.276&0.240\\
150 & 300 &   0.882& 0.892&    0.901 &  0.927&  0.946&  0.973& 0.257& 0.243& 0.222 \\
\hline
\end{tabular}
\label{table2}
\small

\it The columns of loadings and factors report the canonical correlations. PC is the regular principal components method; HWPC represents the heteroskedastic WPC; EWPC uses   $\hSig_u^{-1}$ as the weight matrix. 

\end{center}
\end{table}

We see that the estimation becomes more accurate when we increase the dimensionality.  HWPC improves the regular PC, while the 
 EWPC  gives the best estimation results.

\textbf{Design 2}

Adding a regression term to the model of Design 1, we consider the panel data model with interactive effect: $y_{it}=X_{it}'\beta+\lambda_{i1}f_{1,t}+\lambda_{i,2}f_{2t}+u_{it}$, where the true $\beta=(1,3)'$. The regressors are generated to be dependent on $(f_t,\lambda_i)$:
$$
X_{it,1}=2.5\lambda_{i1}f_{1,t}-0.2\lambda_{i2}f_{2,t}-1+\eta_{it,1}, \quad X_{it,2}=\lambda_{i1}f_{1,t}-2\lambda_{i2}f_{2,t}+1+\eta_{it,2}
$$
where  $\eta_{it,1}$ and $\eta_{it,2}$ are independent  i.i.d. standard normal.

 Both the methods PC  (Bai 2009 and Moon and Weidner  2011)  and   the proposed WPC are carried out to estimate $\beta$ for the comparison. Also compared is the mean squared error of the estimated common components. The simulation is replicated for one hundred times;  results are summarized in Table \ref{table3}. We see that both methods are almost unbiased, while the efficient WPC indeed has significantly smaller standard errors than the regular PC method in the panel model with interactive effects.

  \begin{table}[htdp]
\caption{Method comparison for the panel data with interactive effects, simulation}
\begin{center}
\begin{tabular}{cc|ccccc|ccccc}
 
\hline
&& \multicolumn{5}{c|}{$\beta_1=1$}  &\multicolumn{5}{c}{ $\beta_2=3$}  \\
&& \multicolumn{2}{c}{Mean}  &&\multicolumn{2}{c|}{Normalized SE }& \multicolumn{2}{c}{Mean} & &\multicolumn{2}{c}{Normalized SE} \\
$T$ & $N$ & WPC& PC &&WPC&PC& WPC& PC& &WPC&PC \\
\hline
&&  & & &&&&  &\\

50  & 75  & 1.005 & 1.013 &  & 0.758 & 1.413 & 2.998 & 3.002 &  & 0.744 & 1.472 \\
50  & 100 & 1.005 & 1.010 &  & 0.662 & 1.606 & 2.997 & 2.998 &  & 0.731 & 1.616 \\
50  & 150 & 1.004 & 1.008 &  & 0.964 & 1.913 & 2.999 & 2.999 &  & 0.951 & 1.881 \\
    &     &       &       &  &       &       &       &       &  &       &       \\
100 & 100 & 1.002 & 1.010 &  & 0.550 & 1.418 & 3.000 & 3.003 &  & 0.416 & 1.353 \\
100 & 150 & 1.003 & 1.007 &  & 0.681 & 1.626 & 2.999 & 3.000 &  & 0.611 & 1.683 \\
100 & 200 & 1.002 & 1.005 &  & 0.631 & 1.800 & 3.000 & 3.000 &  & 0.774 & 1.752 \\
    &     &       &       &  &       &       &       &       &  &       &       \\
150 & 100 & 1.003 & 1.006 &  & 0.772 & 1.399 & 3.000 & 2.999 &  & 0.714 & 1.458 \\
150 & 150 & 1.001 & 1.005 &  & 0.359 & 1.318 & 3.000 & 3.001 &  & 0.408 & 1.379 \\
150 & 200 & 1.001 & 1.003 &  & 0.547 & 1.566 & 3.000 & 3.000 &  & 0.602 & 1.762\\

\hline

\end{tabular}
\label{table3}
\small

\it WPC (with weight $\widetilde\Sigma_u^{-1}$) and  PC (existing method) comparison. ``Mean'' is the average of the estimators;   ``Normalized SE'' is   the standard error of the estimators multiplied by  $\sqrt{NT}$.   \end{center}
\end{table}

\section{Empirical Study : Effects of Divorce Law Reforms}\label{law}

This section shows the advantages of our proposed WPC method in a real data application. It demonstrates  the gain of incorporating the estimated $\Sigma_u$ in the panel data estimation and the   efficiency gains  compared to the traditional PC.

\subsection{Real Data Application}

An  important question in sociology is the cause of the sharp increase in the U.S. divorce rate in the  1960s and 1970s. The association between divorce rates and divorce law reforms   has been considered a potential key,   and during 1970s, about three quarters of states in the U.S. liberalized their divorce system, so-called ``no-fault revolution". There is plenty empirical research  regarding the effects of divorce law reforms on the divorce rates (e.g., Peters 1986, Allen 1992), and  statistical significance of these effects has been found (e.g., Friedberg 1998). In other words, states' law reforms are found to have significantly contributed to the increase in state-level divorce rates within the first eight years following reforms. 
 
On the other hand, there is a puzzle about longer effects. Empirical evidence also illustrates the subsequent decrease of the divorce rates  starting from (around) 1975, which is between nine and fourteen  years after the law reforms in most states. So     whether law reforms  continue to contribute to the rate decrease  has been an interesting question.  Wolfers (2006)  studied a  treatment effect panel data model, and    identified    negative effects for the subsequent years. This suggests that, the increase in divorce  following reform and the subsequent decrease  may be two sides of the  same treatment: after earlier dissolution of bad matches  after law reforms, marital relations were gradually affected and changed.  However, it has been argued that Wolfers (2006)'s approach may not capture the complex unobserved heterogeneity. The heterogeneity may exist through an interactive effect, where unobserved common factors may change over time.

 Kim and Oka (2013)  pioneered using interactive effect model for the study: 
\begin{equation}
y_{it}=\sum_{k=1}^KX_{it,k}\beta_k+\lambda_i'f_t+\mu_i+\alpha_t+f(\delta_i, t)+u_{it},
\end{equation}
where  $y_{it}$ is the divorce rate for state $i$ in year $t$; $X_{it, k}$ is a binary regressor, representing the treatment effect   $2k$ years after the reform. Specifically, we observe the law reform year $T_i$ for each state. Then $X_{it,k}=1$ if $2k-1\leq t-T_i\leq 2k$, and zero otherwise. In addition to the interactive effect $\lambda_i'f_t$ as being discussed, the model also contains unobserved state and time effects $(\mu_i, \alpha_t)$ and time trend $f(\delta_i, t)$. For instance, the linear trend defines $f(\delta_i, t)=\delta_i t$ with unknown coefficient $\delta_i.$ Using the regular PC method, Kim and Oka (2013) concluded insignificant $(\beta_5,...,\beta_8)$, that is, the divorce rates after eight years and beyond are not affected by the reforms. However, We argue that  
 using the regular PC method to estimate the model may lose  efficiency   because it ignores the off-diagonal entries. As a result, this can result in wide confidence intervals and possibly conservative conclusions.

We re-estimate Kim and Oka (2013)'s model using the new WPC, and compare with the regular PC.  As a first step, we rewrite the model to fit in the form being considered in this paper.  Introduce the conventional notation: 
$$
\dot{y}_{it}=y_{it}-\frac{1}{T}\sum_{t=1}^Ty_{it}-\frac{1}{N}\sum_{i=1}^Ny_{it}+\frac{1}{NT}\sum_{i=1}^N\sum_{t=1}^Ty_{it}.
$$
Let $\dot{X}_{it,k}, \dot{u}_{it}$ be defined similarly.  
If the time trend $f(\delta_i, t)$ is not present,\footnote{When the time trend is present, we can do a simple projection to eliminate the  time trend, and still estimate the untransformed $\beta$ from the familiar interactive effect model. For instance, suppose $f(\delta_i,t)=\delta_it$. Let 
$
M=(1,2,...,T)'
$ and $P_M=I_T-M(M'M)^{-1}M'$. We  can define 
$
\widetilde{Y}_i=P_M(y_{i1},...,y_{iT})',$ and $ \widetilde{X}_i=P_M(X_{i1},...,X_{iT})',
$
and define $\dot{\widetilde{y}}_{it}$ and $\dot{\widetilde{X}}_{it}$ accordingly from $\widetilde{y}_{it}$ and $\widetilde{X}_{it}$.} under the conventional normalizations $\sum_{i=1}^N\lambda_i=\sum_{t=1}^Tf_t=0,$  $\sum_{i=1}^N\mu_i=\sum_{t=1}^T\alpha_t$, we have
$
\dot{y}_{it}=\dot{X}_{it}'\beta+\lambda_i'f_t+\dot{u}_{it}.
$

The same data as in Wolfers (2006) and Kim and Oka (2013) are used, which contain the divorce rates, state-level reform years and binary regressors from 1956 to 1988 ($T=33$) over  48 states. We fit the models both with and without linear time trend, and apply regular PC and our proposed WPC in each model  to estimate  $\beta$ with confidence intervals.  The number of factors is selected in a data-driven way as in Bai (2009).  His IC and CP  both  suggested ten factors. \footnote{This is the same as in Kim and Oka (2013). We   also tried a few larger values for $r$, and the estimates are  similar, consistent with previous findings that the estimation is robust to over-estimating $r$.} Moreover, for the WPC, the threshold value in the estimated covariance is obtained using the suggested cross-validation procedure in Fan et al. (2013). The estimated $(\beta_1,...,\beta_8)$ and their confidence intervals are summarized in Table \ref{table4}.

Both models produce similar estimates. Interestingly, WPC confirms that the law reforms  significantly contribute to the subsequent decrease of the divorce rates, more specifically,     9-14 years after the reform in the model with linear time trends, and 11-14 years after in the model without linear time trends. In contrast, the regular PC reaches  a more conservative conclusion as it does not capture these significant negative effects. Moreover, both methods show that the effect on the increase of divorce rates for the first 6 years are significant, which is consistent with previous findings in this literature.

 \begin{table}[htdp]
\caption{Method comparison in effects of divorce law reform: real data}
\begin{center}
\begin{tabular}{c|ccccc|c}

\hline
 & \multicolumn{6}{c}{Interactive effect} \\
 &   \multicolumn{2}{c}{WPC}  &  &    \multicolumn{2}{c}{PC} & Relative\\
 & estimate & confidence interval &  & estimate & confidence interval & efficiency \\
 \hline
 &&&&&&\\
First 2 years&   0.014 &     [0.007,    0.021]*  &  &0.018 & [0.0091,    0.028]* & 0.59 \\
3-4 years & 0.034 &    [0.027,    0.041]* &  &  0.042 &    [0.032,    0.053]*  & 0.59 \\
5-6 years &    0.025 &     [0.017,    0.032]* & &0.032 &     [0.022,    0.042]* &  0.58 \\
7-8 years & 0.015 &   [0.007,    0.023]* &   &  0.030 & [0.019,    0.04]*       & 0.56 \\
9-10 years &  -0.006 &    [-0.014,    0.001]    &&0.008 &   [-0.002,    0.018]   & 0.56\\
11-12 years &  -0.008 &    [-0.015,   -0.001]*   &&  0.010 &    [-0.001,    0.02]  & 0.53\\
13-14 years & -0.009 &    [-0.017,   -0.001]*    && 0.005 &    [-0.005,    0.016] &    0.53\\
15 years+ &0.009 &     [0.001,   0.017]*    && 0.031 &     [0.020,    0.042]* &   0.55\\
\hline

 & \multicolumn{6}{c}{Interactive effect+linear trend} \\
 &   \multicolumn{2}{c}{WPC}  &  &    \multicolumn{2}{c}{PC} & Relative\\
 & estimate & confidence interval &  & estimate & confidence interval & efficiency\\
   &&&&&\\

First 2 years &0.014 &  [0.006,    0.021]* & & 0.016 & [0.006,    0.026]*    &  0.55 \\
3-4 years  & 0.032 &     [0.024,    0.039]*  && 0.037 &    [0.026,    0.047]*   & 0.54 \\
5-6 years & 0.018 &     [0.010,    0.026]*   &&0.024 &     [0.012,    0.035]* &  0.54\\
7-8 years  & 0.006 &    [-0.002,    0.014]   & & 0.017 &     [0.005,    0.028]*   & 0.52\\
9-10 years & -0.017 &    [-0.025,   -0.008]*   & & -0.007 &    [-0.019,    0.005]   & 0.52\\
11-12 years & -0.019 &    [-0.028,   -0.010]*   & & -0.006 &    [-0.018,    0.006]  & 0.51\\
13-14 years &-0.021 &    [-0.030,   -0.012]*   && -0.012 &    [-0.025,    0.001]   &  0.50\\
15 years+  & -0.003 &    [-0.012,    0.006]   & & 0.014 &     [0.000,    0.028]*   &0.46\\
\hline

\end{tabular}
\label{table4}
\small

\it 95\% confidence intervals are reported; intervals with * are significant. Relative efficiency is referred to WPC relative to PC,   as estimated $\var(WPC)/\var(PC)$.
\end{center}
\end{table}

We also   report the relative efficiency   using WPC, relative to the regular PC. The reported numbers are   $\var(\text{WPC})/\var(\text{PC})$, where $\var(\text{A})$ calculates the estimated variance of the estimator using method A. It is clear from the table that WPC achieves almost $50\%$ of efficiency gain relative to the regular PC method.

\subsection{Simulated data}

Let us further demonstrate the relative efficiency WPC gains by incorporating the estimated $\Sigma_u^{-1}$ through simulated data. The true parameters are estimated from the real data as described above.  Specifically, we use the first column from Table \ref{table4} (no linear trend) as the true $\beta$, and the corresponding estimated $\Lambda$ as the true loading matrix. We fix $N=48$ as before.  To pertain the actual cross-sectional dependence, in the simulation, the true error terms, factors, and regressors  are  generated as  simple random samples (with replacement) of size $T$ from the estimated residuals, factors and regressors from the real data.   

Simulations are conducted with one hundred replications. The averages and the standard deviations for each estimated component are reported in Table \ref{table5}, representing the bias and standard error.  Also reported is the relative efficiency, defined as $\var(\text{WPC})/\var(\text{PC})$. It is clearly shown in the table that the standard errors of WPC are uniformly smaller than those of PC.  In addition, most of the time WPC also reduces the finite sample bias. The relative efficiency varies from 39\% to 52\%, which illustrates 48\%-61\% efficiency gain.  Overall, after incorporating the error covariance, the performance of the estimator is significantly improved.

 \begin{table}[htdp]
\caption{Method comparison in effects of divorce law reform: simulated data}
\begin{center}
\begin{tabular}{c|cccccccc}

 & \multicolumn{2}{c}{ Bias }    &  &  \multicolumn{2}{c}{Normalized SE} & Relative \\
 & WPC & PC &  & WPC & PC & Efficiency \\
  \hline
 &  \multicolumn{6}{c}{  $T=50$}  \\

First 2 years&   -0.008 &  -0.013 &  &1.077 & 1.714  & 0.393  \\
3-4 years & -0.023  &-0.033  &  &  1.911& 2.694  &  0.494 \\
5-6 years &-0.040  & -0.058  &  &2.743  & 3.821 &  0.525  \\
7-8 years &  -0.054  &  -0.080  &  & 3.429  & 4.899  &  0.501 \\
9-10 years & -0.068 & -0.103 &  &4.017 & 5.633  &  0.501  \\
11-12 years &  -0.073 &  -0.107  &  &4.262  & 6.221 &  0.475  \\
13-14 years&-0.081 &  -0.124  &  &   4.703  &  6.907  & 0.462  \\
15  years+ &-0.090&  -0.133 &  &  5.094  & 7.691 &  0.439  \\
\hline
\\
 &  \multicolumn{6}{c}{  $T=70$}  \\
First 2 years &   -0.002  &  -0.000  &  &   0.927 &  1.449  &  0.408   \\
3-4 years & -0.008 &  -0.008 &  & 1.623&   2.434   & 0.438   \\
5-6 years & -0.021 &   -0.028 &  &    2.434  &  3.420 & 0.505 \\
7-8 years &  -0.030  &   -0.039 &  & 3.246  &   4.579 &   0.507  \\
9-10 years &-0.043  &-0.060    &  &  4.115   & 5.738  &  0.513    \\
11-12 years &   -0.048 & -0.061 &  & 4.579  &  6.492 &  0.501 \\
13-14 years  &  -0.055   & -0.076 &  &   5.101  & 7.245 &  0.495   \\
15  years+ & -0.062  & -0.079  &  & 5.564 &    8.173  &  0.465  \\
\hline

\end{tabular}
\label{table5}

\it ``Normalized SE" is the standard error of the estimator multiplied by $\sqrt{NT}$. The relative efficiency is calculated as the square of the ratio of the third and fourth columns,  estimating $\var(WPC)/\var(PC)$
\end{center}
\end{table}

\section{Conclusion}

The literature on estimating high-dimensional  sparse covariance matrices has targeted on the covariance and inverse covariance directly, and the theoretical results are mostly in an absolute convergence form. We see that the absolute convergence, even though achieving the minimax optimal rate, is often not suitable for statistical inference.  Thus using an estimated high-dimensional covariance matrix as the optimal weight matrix is highly-nontrivial.  We study a new notion of ``weighted convergence" to show that the effect of estimating a high-dimensional covariance matrix is indeed asymptotically negligible.

This paper studies in detail two models that are of increasing importance in applied statistics: approximate factor model and panel data with unobservable interactive effects.  We propose a method of weighted principal components, 
which  uses a high-dimensional weight matrix.  The efficient weight  uses the inverse error covariance matrix. The  EWPC considers both heteroskedasticity and cross-sectional dependence.  It is shown that EWPC uses the optimal weight matrix over the class of WPC estimators thus it is the most efficient.  

  The EWPC is applied to the  year-state divorce rate data. The new method captures the significant (negative) effects from nine to twelve years after the law was reformed, consistent with the previous empirical findings in the social science literature. The estimator is more accurate and produces tighter confidence intervals.
  
  All proofs are given in the supplementary material.

\end{document}